\documentclass[12pt]{amsart}
\usepackage{amssymb,latexsym,amsmath,amsthm,mathrsfs,yfonts}
\usepackage{fullpage,enumerate}
\input xy
\xyoption{all}

\newcommand{\ZZ}{\mathbb Z}

\newcommand{\CC}{\mathbb C}

\newcommand{\NN}{\mathbb N}

\def\TT{\mathbb{T}}

\def\cH{\mathcal H}

\newcommand{\map}{\rightarrow}

\newcommand{\cpt}{\mathbb K}

\def\Aut{\textup{Aut}}

\def\C{\textup{C}}

\def\id{\mathrm{id}}

\def\K{\textup{K}}
\def\KK{\textup{KK}}
\def\RK{\textup{RK}}

\def\sC{\text{$\sigma$-$C^*$}}

\def\lim{{\varprojlim}}
\def\prot{\hat{\otimes}}

\def\dlim{\varinjlim}
\def\ilim{\varprojlim}

\newcommand{\beq}{\begin{eqnarray}}
\newcommand{\beqn}{\begin{eqnarray*}}
\newcommand{\eeq}{\end{eqnarray}}
\newcommand{\eeqn}{\end{eqnarray*}}

\newtheorem{thm}{Theorem}
\newtheorem{lem}{Lemma}

\newtheorem{defn}{Definition}
\newtheorem{rem}{Remark}
\newtheorem{prop}{Proposition}
\newtheorem{ex}{Example}

\begin{document}

\title[Operator algebra quantum homogeneous spaces]
{Operator algebra quantum homogeneous spaces of universal gauge groups}
\author{Snigdhayan Mahanta}
\email{snigdhayan.mahanta@adelaide.edu.au}
\author{Varghese Mathai}
\email{mathai.varghese@adelaide.edu.au}
\address{Department of Pure Mathematics,
University of Adelaide,
Adelaide, SA 5005, 
Australia}

\thanks{{\em Acknowledgements}. Both authors gratefully acknowledge support under the Australian Research Council's {\em Discovery Projects} funding scheme.}



\begin{abstract}
In this paper, we quantize universal gauge groups such as $SU(\infty)$, as well as their homogeneous spaces, in the $\sC$-algebra setting. More precisely, we propose concise definitions of $\sC$-quantum groups and $\sC$-quantum homogeneous spaces and explain these concepts here. At the same time, we put these definitions in the mathematical context of countably compactly generated spaces as well as $C^*$-compact quantum groups and homogeneous spaces.  
We also study the representable 
K-theory of these spaces and compute it for the quantum homogeneous spaces associated to the universal gauge group $SU(\infty)$.

\end{abstract}

\maketitle


If $H$ is a compact and Hausdorff topological group, then the $C^*$-algebra of all continuous functions $\C(H)$ admits a comultiplication map $\Delta: \C(H)\map \C(H)\prot\C(H)$ arising from the multiplication in $H$. This observation motivated Woronowicz (see, for instance, \cite{W}), amongst others such as Soibelman \cite{S}, to introduce the notion of a {\it $C^*$-compact quantum group} in the setting of operator algebras as a unital $C^*$-algebra with a coassociative comultiplication, satisfying a few other conditions.  If the group $H$ is only locally compact then the situation becomes significantly more difficult. One of the reasons is that the multiplication map $m:H\times H\map H$ is no longer a proper map and one needs to introduce multiplier algebras of $C^*$-algebras to obtain a comultiplication, see for instance, Kustermans-Vaes \cite{KV}. For an excellent and thorough introduction to this theory the readers are referred to, for instance, \cite{KVVVDW}. In the sequel we show that if $H=\dlim_n H_n$ is a countably compactly generated group, i.e., if $H_n\subset H_{n+1}$ are compact and Hausdorff topological groups for all $n\in\NN$ and if $H$ is the direct limit, then a story similar to the compact group case goes through using the general framework of $\sC$-algebras as systematically developed by Phillips \cite{P1,P2}, motivated by some earlier work by Arveson, Mallios, Voiculescu, amongst others. There is a clean formulation of, what we call, {\it $\sC$-quantum groups}, which are noncommutative generalizations of $\C(H)$. Examples of countably compactly generated groups are $U(\infty)=\dlim_n U(n)$, $SU(\infty) =\dlim_n SU(n)$, where $U(n)$ (resp. $SU(n)$) are the unitary (resp. special unitary) groups. They are also known in the physics literature as {\em universal gauge groups}, see Harvey-Moore \cite{HM} and Carey-Mickelsson \cite{CM}. Such spaces are not locally compact and hence the existing literature on quantum groups cannot handle them. Moreover, locally compact groups that are not compact, are also not countably compactly generated. We also discuss in detail the interesting example of the quantum version of the universal special unitary group, $C(SU_q(\infty))$. 

A {\it pro $C^*$-algebra} is an inverse limit of $C^*$-algebras and $*$-homomorphisms, where the inverse limit is constructed inside the category of all topological $*$-algebras and continuous $*$-homomorphisms. For the general theory of topological $*$-algebras one may refer to, for instance, \cite{M}. The underlying topological $*$-algebra of a pro $C^*$-algebra is necessarily complete and Hausdorff. It is not a $C^*$-algebra in general; it would be so if, for instance, the directed set is finite. If the directed set is countable, then the inverse limit is called a {\it $\sC$-algebra}. One can choose a linearly directed cofinal subset inside any countable directed set and the passage to a cofinal subsystem does not change the inverse limit. Therefore, we shall always identify a $\sC$-algebra $A\cong\ilim_n A_n$, where $n\in\NN$. The inverse limit could have also been constructed inside the category of $C^*$-algebras; however, the two results will not agree. For instance, if $H=\dlim_n H_n$ as above, then the inverse limit $\ilim_n \C(H_n)$ inside the category of topological $*$-algebras is $\C(H)$, whereas that inside the category of $C^*$-algebras is $\C_b(H)$, i.e., the norm bounded functions on $H$. It is known that $\C_b(H)\cong\C(\beta H)$, where $\beta H$ is the Stone--{\v{C}}ech compactification of $H$. Therefore, if one wants to model a space via its algebra of all continuous functions then the former inverse limit is the appropriate one. {\em Henceforth, the inverse limits are always constructed inside the category of topological $*$-algebras}. It is known that any $*$-homomorphism between two pro $C^*$-algebras is automatically continuous, provided the domain is a $\sC$-algebra (see Theorem 5.2. of \cite{P1}). Furthermore, the category of commutative and unital $\sC$-algebras with unital $*$-homomorphisms (automatically continuous) is contravariantly equivalent to the category of countably compactly generated and Hausdorff spaces with continuous maps via the functor $X\mapsto \C(X)$ (see Proposition 5.7. of \cite{P1}). If $A\cong\ilim_n A_n$, $B\cong\ilim_n B_n$ are two $\sC$-algebras, then the minimal tensor product is defined to be $A\prot_{\textup{min}} B = \ilim_n A_n\prot_{\textup{min}} B_n$. {\em Henceforth, $A\prot B$ will always denote the minimal or spatial tensor product between $\sC$-algebras}.

We next outline the contents of the paper. \S\ref{sect:qgroup} initiates the concept of a $\sC$-quantum group, where  the interesting example of the quantum version of the universal special unitary group, $C(SU_q(\infty))$, is discussed in detail. \S\ref{sect:qhom}  initiates the concept of a $\sC$-quantum homogeneous space, where some interesting examples of quantum version of 
homogeneous spaces associated to the universal special unitary group, $SU(\infty)$, are discussed in detail. \S\ref{sect:RK} contains the computation of the representable K-theory of $C(SU_q(\infty))$ as well as some of
the quantum homogeneous spaces associated to it.

\section{$\sC$-quantum groups}\label{sect:qgroup}
In  this section, we define the concept of a $\sC$-quantum group and explain it here. We also discuss in detail the interesting example of the quantum version of the universal special unitary group, $C(SU_q(\infty))$.

If $H$ is a countably compactly generated and Hausdorff topological group, although the multiplication map $m:H\times H\map H$ is not proper, we get an induced comultiplication map $m^*:\C(H)\map\C(H\times H)\cong\C(H)\prot\C(H)$, which will be coassociative owing to the associativity of $m$. Motivated by the definition of Woronowicz (see also Definition 1 of \cite{KV}), we propose:

\begin{defn}
A unital $\sC$-algebra $A$ is called a $\sC$-quantum group if there is a unital $*$-homomorphism $\Delta: A\map A\prot A$ which satisfies coassociativity, i.e., $(\Delta\prot\id)\Delta = (\id\prot\Delta)\Delta$ and such that the linear spaces $\Delta(A)(A\prot 1)$ and $\Delta(A)(1\prot A)$ are dense in $A\prot A$.
\end{defn}

\begin{lem} 
Let $\{A_n,\theta_n: A_n\map A_{n-1}\}_{n\in\NN}$ be a countable inverse system of $C^*$-algebras and let $B_n\subset A_n$ be dense subsets for all $n$ such that $\theta_n(B_n)\subset B_{n-1}$. Then $\ilim_n B_n$ is a dense subset of the $\sC$-algebra $\ilim_n A_n$.
\end{lem}

\begin{proof}
The assertion follows from the Corollary to Proposition 9 in \S 4-4 of \cite{Bou}.
\end{proof}

\begin{ex} \label{example}
Obviously, any $C^*$-compact quantum group is a $\sC$-quantum group. Let $\{A_n,\theta_n: A_n\map A_{n-1}\}_{n\in\NN}$ be a countable inverse system of $C^*$-compact quantum groups with $\theta_n$ surjective and unital for all $n$. Furthermore, let us assume that the comultiplication homomorphisms $\Delta_n$ form a morphism of inverse systems of $C^*$-algebras $\{\Delta_n\}:\{A_n\}\map\{A_n\prot A_n\}$. Then $(A,\Delta)=(\ilim_n A_n,\ilim_n\Delta_n)$ is a $\sC$-quantum group. Indeed, the density of the linear spaces $\Delta(A)(A\prot 1)$ and $\Delta(A)(1\prot A)$ inside $A\prot A$ follow from the above Lemma.
\end{ex}

Our next goal is to outline the construction of the quantum universal special unitary group, $C(SU_q(\infty))$.
Recall that for $q\in (0,1)$, the $C^*$-algebra $C(SU_{q}(n))$ is the universal $C^*$-algebra 
generated by $n^{2}+2$ elements $G_n:=\{u^n_{ij}: \, i, j =1, \ldots , n\}\cup\{0,1\}$, which satisfy the following
relations

\begin{equation}\label{relations0}
0^*=0^2=0,\;\;1^*=1^2=1, \;\;01=0=10, \;\;1u^n_{ij} = u^n_{ij}1 = u^n_{ij},\;\; 0u^n_{ij}=u^n_{ij}0=0 \text{  for all $i,j$}
\end{equation}

\begin{equation}\label{relations1}
 \sum_{k=1}^{n}u^n_{ik}(u^{n}_{jk})^*=\delta_{ij}1, \qquad \sum_{k=1}^{n}(u^{n}_{ki})^*u^n_{kj}
=\delta_{ij}1
\end{equation}
\begin{equation}\label{relations2}
\sum_{i_{1}=1}^{n}\sum_{i_{2}=1}^{n} \cdots \sum_{i_{n}=1}^{n}
E_{i_{1}i_{2}\cdots i_{n}}u^n_{j_{1}i_{1}}\cdots u^n_{j_{n}i_{n}} =
E_{j_{1}j_{2}\cdots j_{n}}1
\end{equation}
where
\begin{displaymath}
\begin{array}{lll}
E_{i_{1}i_{2}\cdots i_{n}}&:=& \left\{\begin{array}{lll}
                                 0 &\text{whenever}& i_{1},i_{2},\cdots, i_{n} \text{ are not
distinct;} \\
                                 (-q)^{\ell(i_{1},i_{2},\cdots,i_{n})}.
                                 \end{array} \right. 
\end{array}
\end{displaymath}
Here $\delta_{ij}1=0$ if $i\neq j$, where $0$ denotes the zero element in the generating set of $C(SU_q(n))$ and, for any permutation $\sigma$, $\ell(\sigma)=\textup{Card}\{(i,j)\,|\, i<j,\sigma(i)>\sigma(j)\}$. 
The $C^{*}$-algebra $C(SU_{q}(n))$ has a $C^*$-compact quantum group
structure with the comultiplication $\Delta$ given by 
\begin{displaymath}
\Delta(0):=0\otimes 0,\;\; \Delta(1):=1\otimes 1 \text{  and   } \Delta(u^n_{ij}):= \sum_{k} u^n_{ik}\otimes u^n_{kj}.
\end{displaymath} It is known that $C(SU_q(n))$ is a type-$\textup{I}$ $C^*$-algebra \cite{Br}, whence it is nuclear. Therefore, there is a unique choice for the $C^*$-tensor product in the definition of the comultiplication. There is a surjective $*$-homomorphism $\theta_n:C(SU_{q}(n)) \to C(SU_{q}(n-1))$ defined on the generators by
 \begin{eqnarray*}
 \theta_n(x) &:=& x \text{  if $x=0,1$}\\
 \theta_n(u^n_{ij}) &:=& \left\{\begin{array}{ll}
                        u^{n-1}_{ij}&\text{if}~ 1\leq i,j \leq n-1 ,\\
                        \delta_{ij}1 & \text{otherwise},
                        \end{array} \right.
\end{eqnarray*} such that the following diagram commutes for all $n\geqslant 2$

\beqn \label{comult}
\xymatrix{
C(SU_q(n))\ar[rrr]^{\Delta_n} \ar[d]_{\theta_n} &&& C(SU_q(n))\prot C(SU_q(n)) \ar[d]^{\theta_n\prot\theta_n}\\
C(SU_q(n-1))\ar[rrr]^{\Delta_{n-1}} &&& C(SU_q(n-1))\prot C(SU_q(n-1)).
}
\eeqn One can verify this assertion by a routine calculation on the generators. Consequently, for $n\geqslant 2$ the families $\{C(SU_{q}(n)), \theta_n\}$ and $\{C(SU_q(n))\prot C(SU_q(n)),\theta_n\prot\theta_n\}$ form countable inverse systems of $C^*$-algebras and $\{\Delta_n\}:\{C(SU_q(n))\}{\map} \{C(SU_q(n))\prot C(SU_q(n))\}$ becomes a morphism of inverse systems of $C^*$-algebras. We construct the underlying $\sC$-algebra of the universal quantum gauge group as the inverse limit $$C(SU_q(\infty))=\ilim_n C(SU_{q}(n)).$$ In fact, $C(SU_{q}(\infty))$ is a $\sC$-quantum group, since it is the inverse limit of $C^*$-compact quantum groups, where the comultiplication $\Delta$ on $C(SU_{q}(\infty))$ is defined to be $\Delta =\ilim_n \Delta_n$ (see the Example above).

%
%
%
%
%
%

If $G$ is a set of generators and $R$ a set of relations, such that the pair $(G,R)$ is {\it admissible} (see Definition 1.1. of \cite{B}), then one can always construct a universal $C^*$-algebra $C^*(G,R)$. For instance, the universal $C^*$-algebra generated by the set $\{1,x\}$, subject to the relations $\{1^*=1^2=1$, $1x=x1=x$, $x^* x= 1=xx^*\}$, is isomorphic to $\C(S^1)$. The generators and relations of $C(SU_q(n))$ described above are also admissible. 

\begin{rem}
All matrix $C^*$-compact quantum groups considered, for instance, in \cite{W1,W}, such that the relations put a bound on the norm of each generator, are of the form $C^*(G,R)$, where $(G,R)$ is an admissible pair of generators and relations.
\end{rem}

Let $\{(G_i, R_i)\}_{i\in\NN}$ be a countable family of admissible pairs of generators and relations, so that $C^*(G_i, R_i)$ exist for all $i$. Let $F(G)$ denote the associative nonunital complex $*$-algebra (freely) generated by the concatenation of the elements of $G\coprod G^*$ and finite $\CC$-linear combinations thereof, where $\coprod$ denotes disjoint union and $G^*=\{g^*\,|\, g\in G\}$ (formal adjoints). We call a relation in $R$ {\it algebraic} if it is of the form $f=0$ (or can be brought to that form), where $f\in F(G)$. For instance, if $G=\{1,x\}$, then $x^*x=1$ is algebraic, whereas $\|x\|\leqslant 1$ is not. If $(G,R)$ is a pair of generators and relations, then a {\it representation} $\rho$ of $(G,R)$ in a (pro) $C^*$-algebra $B$ is a set map $\rho:G\map B$, such that $\rho(G)$ satisfies the relations $R$ inside $B$. If $(G,R)$ is a {\it weakly admissible} pair of generators and relations (see Definition 1.3.4. of \cite{P2}), then one can construct the universal pro $C^*$-algebra $C^*(G,R)$ (see Proposition 1.3.6. of ibid.). It is known that any combination (even the empty set) of algebraic relations is weakly admissible (see Example 1.3.5.(1) of ibid.).

\noindent
We further make the following hypotheses:

\begin{enumerate}[(a)]
\item \label{Hypa} There are surjective maps $\theta_i: G_i\map G_{i-1}$, so that one may form the inverse limit in the category of sets $G=\ilim_i G_i$, with canonical projection maps $p_i:G\map G_i$. We also require the surjections $\theta_i$ to admit sections $s_{i-1}:G_{i-1}\map G_i$ satisfying $\theta_i\circ s_{i-1}=\id_{G_{i-1}}$, so that we get canonical splittings $\gamma_i: G_i\map G$ satisfying $p_i\circ\gamma_i =\id_{G_i}$. The map $\gamma_i$ sends $g_i\map\{h_j\}$, where

\beq \label{split}
h_j = \begin{cases}
g_i \text{  if $j=i$,}\\
\theta_{i-n+1}\circ\cdots\circ \theta_{i} (g_i) \text{  if $j=i-n$, $n>0$,}\\
s_{i+m-1}\circ\cdots\circ s_i (g_i) \text{  if $j=i+m$, $m>0$.}
\end{cases}
\eeq


\item \label{Hypb} We require that for all $i$ the iterated applications of $\theta_j$'s  and $s_k$'s on $G_i$ satisfy $R_i$ for all $j\leqslant i$ and $k\geqslant i$. 



\end{enumerate} The surjective maps $\theta_i$ induce surjective $*$-homomorphisms $\theta_i:C^*(G_i,R_i)\map C^*(G_{i-1},R_{i-1})$; consequently, $\{C^*(G_i,R_i),\theta_i\}_{i\in\NN}$ forms a countable inverse system of $C^*$-algebras. We may form the inverse limit $\ilim_i C^*(G_i,R_i)$, which is by construction a $\sC$-algebra. Let $(G,R)$ be a pair of generators and relations, where $G=\ilim_i G_i$ and $R$ denotes the set of relations $\{\gamma_i(G_i) \text{ satisfies } R_i \text{ for all $i$}\}$. A {\it representation} $\rho$ of $(G,R)$ in a (pro) $C^*$-algebra $B$ is a set map $\rho:G\map B$, such that $\rho\circ\gamma_i(G_i)$ satisfies $R_i$ inside $B$ for all $i$. We assume that $(G,R)$ is a weakly admissible pair, so that one can construct the universal pro $C^*$-algebra $C^*(G,R)$.

\begin{thm}\label{thm:inverselim}
There is an isomorphism of pro $C^*$-algebras $C^*(G,R)\cong \ilim_i C^*(G_i,R_i)$.
\end{thm}

\begin{proof}
It suffices to show that $\ilim_i C^*(G_i,R_i)$ is a universal representation of $(G,R)$, i.e., there is a map $\iota:G\map\ilim_i C^*(G_i,R_i)$ such that $\iota\circ\gamma_i(G_i)$ satisfies $R_i$ inside $\ilim_i C^*(G_i,R_i)$ for all $i$ and given any representation $\rho$ of the pair $(G,R)$ in a pro $C^*$-algebra $B$, there is a unique continuous $*$-homomorphism $\kappa: \ilim_i C^*(G_i,R_i)\map B$ making the following diagram commute:

\beqn
\xymatrix{
G=\ilim_i G_i\ar[rrd]_{\rho}\ar^{\iota}[rr] &&\ilim_i C^*(G_i,R_i) \ar[d]^{\kappa}\\
&& B.
}
\eeqn 
The map $\iota: G\map\ilim_i C^*(G_i,R_i)$ is defined as $g\mapsto \{p_i(g)\}$, which is a representation of $(G,R)$ due to the Hypothesis \eqref{Hypb} above. The construction of the universal pro $C^*$-algebra $C^*(G,R)$ (resp. $C^*$-algebra $C^*(G_i,R_i)$) is defined via a certain Hausdorff completion of $F(G)$ (resp. $F(G_i)$) with respect to representations in pro $C^*$-algebras (resp. $C^*$-algebras) satisfying $R$ (resp. $R_i$). The surjective maps $\theta_i$ induce $*$-homomorphisms $\theta_i: F(G_i)\map F(G_{i-1})$, whence we may construct  the $*$-algebra $\ilim_i F(G_i)$ (purely algebraic inverse limit). By the above Lemma it suffices to define $\kappa$ on coherent sequences of the form $\{w_i\}\in\ilim_i F(G_i)$, which then extends uniquely to a $*$-homomorphism on the entire $\ilim_i C^*(G_i,R_i)$. Thanks to the maps $\rho\circ\gamma_i:G_i\map B$, $\rho$ extends uniquely to a $*$-homomorphism $\ilim_i F(G_i)\map B$. Now there is a unique choice for $\kappa(\{w_i\})$ forced by the compatibility requirement, i.e., $\kappa(\{w_i\})=\rho(\{w_i\})$. By construction $\kappa$ is a $*$-homomorphism and it is automatically continuous, since $\ilim_i C^*(G_i,R_i)$ is a $\sC$-algebra. 
\end{proof}

In the example of $C(SU_q(\infty))$, one could try to define the section maps $s_{n-1}: G_{n-1}\map G_n$ as $$ 0\mapsto 0,\;\; 1\mapsto 1,\;\; u^{n-1}_{ij}\mapsto u^n_{ij}.$$ But the Hypothesis \eqref{Hypb} will not be satisfied and hence the above Theorem is unfortunately not applicable. However, the Theorem could be of independent interest as it can be applied to inverse systems, where the structure $*$-homomorphisms admit sections (also $*$-homomorphisms).

\noindent
Let $G_n:=\{w^n_{ij}: \, i, j =1, \ldots , n\}\cup\{0,1\}$ be a set of generators satisfying the relations $R_n$ $$0^*=0^2=0,\;\;1^*=1^2=1, \;\;01=0=10, \;\;1w^n_{ij} = w^n_{ij}1 = w^n_{ij},\;\; 0w^n_{ij}=w^n_{ij}0=0,\;\, \|w^n_{ij}\|\leqslant 1$$ for all $i,j$. The pair $(G_n,R_n)$ is an admissible pair for all $n$, so that there is a universal $C^*$-algebra $C^*(G_n,R_n)$. There are surjective maps $\theta_n:G_n\map G_{n-1}$ given by 

 \begin{eqnarray*}
 \theta_n(x) &:=& x \text{  if $x=0,1$}\\
 \theta_n(w^n_{ij}) &:=& \left\{\begin{array}{ll}
                        w^{n-1}_{ij}&\text{if}~ 1\leq i,j \leq n-1 ,\\
                        \delta_{ij}1 & \text{otherwise}.
                        \end{array} \right.
\end{eqnarray*} making $\{C^*(G_n,R_n),\theta_n\}$ an inverse system of $C^*$-algebras and surjective $*$-homomorphisms. There are obvious sections $s_{n-1}: G_{n-1}\map G_n$ sending $0\mapsto 0$, $1\mapsto 1$ and $w^{n-1}_{ij}\mapsto w^n_{ij}$ giving rise to maps $\gamma_n:G_n\map G=\ilim_n G_n$ as described above (see Equation \eqref{split}). There are surjective $*$-homomorphisms $\pi_n: C^*(G_n,R_n)\map C(SU_q(n))$ for all $n\geqslant 2$ given on the generators by $\pi_n(x)=x$ for $x=0,1$ and $\pi_n(w^n_{ij})=u^n_{ij}$, which produce a morphism of inverse systems $\{\pi_n\}:\{C^*(G_n,R_n)\}\map\{C(SU_q(n))\}$. Indeed, it follows from the Relations \eqref{relations0}, \eqref{relations1} and \eqref{relations2} that the norms of the generators of $C(SU_q(n))$ do not exceed $1$ in any representation. Consequently, there is a surjective $*$-homomorphism of $\sC$-algebras (see \cite{P4} 1.6. Lemma) $$\ilim_n\pi_n:\ilim_n C^*(G_n,R_n)\map C(SU_q(\infty)).$$
However, the authors cannot provide a good description of the kernel at the moment. Let us set $G=\ilim_n G_n$ and let $R$ denote the set of relations $\{\gamma_n(G_n) \text{ satisfies $R_n$ for all $n$}\}.$ Note that $\|x\|\leqslant 1$ viewed as a relation for a representation in a pro $C^*$-algebra $B$ means that $p(x)\leqslant 1$ for all $C^*$-seminorms $p$ on $B$. The family of pairs $(G_n,R_n)$ satisfy Hypotheses \eqref{Hypb} and the pair $(G,R)$ is weakly admissible (see Example 1.3.5.(2) of \cite{P2}), so that the above Theorem applies, i.e., $\ilim_n C^*(G_n, R_n)\cong C^*(G,R)$. As a corollary, we deduce that the elements of $(\ilim_n\pi_n)(G)$ provide explicit generators of $C(SU_q(\infty))$.



\section{$\sC$-quantum homogeneous spaces}\label{sect:qhom}
In this section we define the concept of a $\sC$-quantum homogeneous space and explain it here. We also discuss in detail the interesting examples of the quantum versions of the homogeneous spaces associated to the universal special unitary group, $SU(\infty)$.

Let $\{A_n,\theta_n: A_n\map A_{n-1}\}_{n\in\NN}$ and $\{B_n,\psi_n: B_n\map B_{n-1}\}_{n\in\NN}$ be countable inverse systems of $C^*$-compact quantum groups with $\theta_n$ and $\psi_n$ surjective and unital for all $n$. Furthermore, let us assume that the comultiplication homomorphisms $\Delta^A_n, \Delta^B_n$ form morphisms of inverse systems of $C^*$-algebras $\{\Delta^A_n\}:\{A_n\}\map\{A_n\prot A_n\}$ and $\{\Delta^B_n\}:\{B_n\}\map\{B_n\prot B_n\}$. Then $(A,\Delta^A)=(\ilim_n A_n,\ilim_n\Delta^A_n)$ and $(B,\Delta^B)=(\ilim_n B_n,\ilim_n\Delta^B_n)$ are $\sC$-quantum groups by the discussion in the above section (see Example \ref{example}).



Suppose now that there are compatible $*$-homomorphisms, $\theta'_n : A_n \to B_n$, that is, such that the following 
diagrams commute
\beq\label{diag:comm2}
\xymatrix{
 A_n\ar[rr]^{{\theta'}_n} \ar[d]_{\theta_n} && B_{n} \ar[d]^{\psi_{n}}\\
A_{n-1}\ar[rr]^{{\theta'}_{n-1}} && B_{n-1}.
}
\eeq
and 
\beq\label{diag:comm3}
\xymatrix{
 A_n\ar[rr]^{\Delta^A_n} \ar[d]_{\theta'_n} && A_n \prot A_{n} \ar[d]^{\theta'_{n}\prot\theta'_n}\\
B_{n}\ar[rr]^{\Delta^B_{n}} && B_{n} \prot B_{n}.
}
\eeq

Then, after Nagy \cite{N} one calls the $C^*$-subalgebras
$$
\cH_n = \left\{f\in A_n\big| (\theta_n' \prot \id)\Delta^A_n(f) = 1\prot f \right\}\subset A_n
$$ 
as {\it $C^*$-compact quantum homogeneous spaces}  for all $n\in\NN$. It is pointed out by Nagy that a parallel theory can be developed for 
$
\tilde\cH_n = \left\{f\in A_n\big| ( \id \prot \theta_n' )\Delta^A_n(f) = f \prot 1 \right\}\subset A_n
$. By assumption, one has the following commutative diagram for all $n$,
\beq\label{diag:comm}
\xymatrix{
 A_n\ar[rr]^{\Delta_n} \ar[d]_{\theta_n} && A_n \prot A_{n} \ar[d]^{\theta_{n}\prot\theta_n}\\
A_{n-1}\ar[rr]^{\Delta_{n-1}} && A_{n-1} \prot A_{n-1}.
}
\eeq
A similar commutative diagram holds for $\{B_n,\psi_n: B_n\map B_{n-1}\}_{n\in\NN}$. Then we have 

\begin{lem} In the notation above, 
the $*$-homomorphism $\theta_n: A_n \to A_{n-1}$ restricts to a $*$-homomorphism of $C^*$-compact quantum homogeneous spaces $\cH_n\map \cH_{n-1}$ for all $n\in\NN$.
\end{lem}

\begin{proof}
Let $f\in \cH_n$. Then $(\theta_n' \prot \id)\Delta_n^A(f) = 1\prot f$ and we compute,
\begin{align*}
(\theta_{n-1}' \prot \id)\Delta_{n-1}^A(\theta_n(f)) & = (\theta_{n-1}' \prot \id)(\theta_{n}\prot\theta_n)\Delta_n^A(f)  \\
& =  (\theta_{n-1}' \prot \id)(\id \prot \theta_n)(\theta_{n}\prot \id) \Delta_n^A(f) \\
& =  (\theta_{n-1}' \prot \id)(\id \prot \theta_n) 1 \prot f \\
& =  1 \prot \theta_n(f),
\end{align*}
showing that $\theta_n(f) \in \cH_{n-1}$. 
\end{proof} 

\noindent
Now we define
\beq
\cH :=\ilim_n \cH_n
\eeq
to be a $\sC$-{\em quantum homogeneous space}, where the inverse limit is once again taken inside the category of topological $*$-algebras. 

\begin{rem}We remark here that this definition can be generalized as follows: Let $(A, \Delta^A)$
and $(B, \Delta^B)$ be $\sC$-quantum groups, and $\theta':A\longrightarrow B$ be a $*$-homomorphism.
Then one can also call the $\sC$-subalgebra
$$
\cH = \left\{f\in A\big| (\theta' \prot \id)\Delta^A(f) = 1\prot f \right\}\subset A
$$ 
a { $\sC$-quantum homogeneous space}. However we will not be discussing these here.
\end{rem}

Our next goal is to outline the construction of the quantum homogeneous space associated to the universal gauge group $SU(\infty)$. Recall that there is a surjective $*$-homomorphism $\theta_n:C(SU_{q}(n)) \to C(SU_{q}(n-1))$ defined on the generators by
 \begin{eqnarray*}
 \theta_n(x) &:=& x \text{  if $x=0,1$}\\
 \theta_n(u^n_{ij}) &:=& \left\{\begin{array}{ll}
                        u^{n-1}_{ij}&\text{if}~ 1\leq i,j \leq n-1 ,\\
                        \delta_{ij}1 & \text{otherwise}.
                        \end{array} \right.
\end{eqnarray*}

\noindent
Then the quantum spheres (cf. \cite{N}) are by fiat 
$$C(S^{2n-1}_q) = \left\{f\in C(SU_{q}(n))\big| (\theta_n \prot \id)\Delta_n(f) = 1\prot f \right\},$$
and come with induced $*$-homomorphisms
$$
\theta_n : C(S^{2n-1}_q) \to C(S^{2n-3}_q).
$$
Then the quantum homogeneous space associated to the universal gauge group 
$SU(\infty)$ is defined to be
$$
C(S^\infty_q) := \ilim_n C(S^{2n-1}_q).
$$
It is shown in \cite{Sheu2} that $C(S^{2n-1}_q)$ is isomorphic to a groupoid $C^*$-algebra, which is independent of $q$.

\begin{ex}
Another example is that of the $C^*$-quantum projective space (see, for instance, \cite{Sheu3}),
$$
C(\CC P^{n}_q) = C^*(\{v_i^*v_j\,|\,v_i=u^{n+1}_{(n+1)i},\,  v_j=u^{n+1}_{(n+1)j},\, 1\leqslant i,j\leqslant n+1\})\subset C(S_q^{2n+1}).
$$

\noindent
Moreover, there is a short exact sequence of $C^*$-algebras relating the $C^*$-quantum projective spaces (see Corollary 2 of \cite{Sheu3}), viz.,

\beq \label{ProjDecomp}
0\map \cpt \map C(\CC P^n_q)\map C(\CC P^{n-1}_q)\map 0,
\eeq for $n\geqslant 1$ and $C(\CC P^0_q)\simeq \CC$. We define the $C^*$-quantum infinite projective space as 
$$
C(\CC P^{\infty}_q) = \ilim_n C(\CC P^{n}_q).
$$

\end{ex}

\section{Representable $\K$-theory of $\sC$-quantum homogeneous spaces}\label{sect:RK}

The appropriate $\K$-theory for $\sC$-algebras is representable $\K$-theory, denoted by $\RK$.
In this section, we compute the representable K-theory of $C(SU_q(\infty))$ as well as some of
the quantum homogeneous spaces associated to it.

 The $\RK$-theory agrees with the usual $\K$-theory of $C^*$-algebras if the input is a $C^*$-algebra and many of the nice properties that $\K$-theory satisfies generalise to $\RK$-theory. Let us briefly recall some of the basic facts about $\sC$-algebras and $\RK$-theory after Phillips \cite{P3} and Weidner \cite{Weid}.

\begin{enumerate}
\item The $\RK$-theory is homotopy invariant and satisfies Bott $2$-periodicity.

\item \label{K} If $A$ is a $C^*$-algebra, then there is a natural isomorphism $\RK_i(A)\cong\K_i(A)$. 

\item There is a natural isomorphism $\RK_i(A\prot \cpt)\cong \RK_i(A)$, where $\cpt$ denotes the algebra of compact operators on a separable Hilbert space.

\item \label{Milnor} If $\{A_n\}_{n\in\NN}$ is a countable inverse system of $\sC$-algebras with surjective homomorphisms (which can always be arranged), then the inverse limit exists as a $\sC$-algebra and there is a Milnor $\lim^1$-sequence

\beqn
0\map \lim^1_n \RK_{1-i}(A_n)\map \RK_i(\lim_n  A_n)\map \lim_n \RK_i(A_n)\map 0.
\eeqn
\end{enumerate}

\noindent
Here we recall that Sheu \cite{Sheu} and Soibelman--Vaksman \cite{S-V} have computed the $\K$-theory
of the $C^*$-quantum spheres, viz., 

$$\K_0(C(S^{2n-1}_q)) \simeq\ZZ \text{    and    } \K_1(C(S^{2n-1}_q)) \simeq \ZZ.$$

\begin{thm}
$\RK_0(C(S^\infty_q))\simeq \ZZ$ and  $\RK_1(C(S^\infty_q))\simeq \{0\}$.
\end{thm}

\begin{proof}
There is a short exact sequence for all $n$ (see Corollary 8 of \cite{Sheu}),
$$
0\to C(\TT)\prot \cpt \longrightarrow C(S^{2n-1}_q) \stackrel{\theta_n}{\longrightarrow} C(S^{2n-3}_q) \to 0.
$$ This gives rise to a $6$-term exact sequence involving the topological $\K$-theory groups

\beq
\xymatrix{
\ZZ\ar[r]^{d_1}&\ZZ\ar[r]^{d_2} &\ZZ\ar[d]^{d_3}\\
\ZZ\ar[u]^{d_6}&\ZZ\ar[l]^{d_5} &\ZZ\ar[l]^{d_4}
}
\eeq Sheu argues that $d_1=0$ (see Section 7 of ibid.) from which it follows that $d_3=d_5=0$ and $d_2=d_4=d_6=\id$ or $-\id$. The differential $d_2$ (resp. $d_5$) is the homomorphism induced by $\theta_n$ between the $\K_0$-groups (resp. $\K_1$-groups). By properties \eqref{K} and \eqref{Milnor} above, one obtains the following exact sequence of abelian groups

\beqn
0\map \lim^1_n \K_{1-i}(C(S^{2n-1}_q))\map \RK_i(C(S^\infty_q))\map \lim_n \K_i(C(S^{2n-1}_q))\map 0.
\eeqn Now $\ilim_n\K_0(C(S^{2n-1}_q))\simeq \ZZ$ since all the connecting homomorphisms are isomorphisms and the $\lim^1$-term vanishes as the connecting homomorphisms between the $\K_1$-groups are all zero, whence $\RK_0(C(S^\infty_q))\simeq\ZZ$. 
Similarly, $\ilim_n\K_1(C(S^{2n-1}_q))\simeq\{0\}$ and the $\lim^1$-term involving the $\K_0$-groups vanishes as the Mittag--Leffler condition is satisfied. It follows that $\RK_1(C(S^\infty_q))\simeq \{0\}$.

\end{proof}


\noindent
Let us now compute the $\RK$-theory of $C(\CC P^\infty_q))$. The following result is presumably well-known, cf. 
\cite{Sheu3}.

\begin{prop}
$\K_0(C(\CC P^n_q))\simeq \ZZ^{n+1}$ and $\K_1(C(\CC P^n_q))\simeq \{0\}$.
\end{prop}

\begin{proof}
We argue by induction on $n$. For $n=0$ the assertion is true since $C(\CC P^0_q)\simeq \CC$. Let us set $A_n = \K_0(C(\CC P^n_q))$ and $B_n = \K_1(C(\CC P^n_q))$, so that $A_0=\ZZ$ and $B_0=\{0\}$. The $6$-term sequence associated to the short exact sequence \eqref{ProjDecomp} gives

\beq
\xymatrix{
\ZZ\ar[r]& A_n\ar[r] & A_{n-1}\ar[d]\\
B_{n-1}\ar[u]& B_n\ar[l] & 0\ar[l]
}
\eeq By the induction hypothesis we obtain $A_{n-1} =\ZZ^n$ and $B_{n-1}=\{0\}$. It follows immediately that $B_n=\{0\}$ and $A_n$ fits into a short exact sequence $$0\map\ZZ\map A_n\map\ZZ^n\map 0.$$ Since $\ZZ^n$ is a projective $\ZZ$-module, this sequence splits, whence $A_n=\ZZ^{n+1}$. 

\end{proof}

\begin{thm}
$\RK_0(C(\CC P^\infty_q))\simeq \ilim_n \ZZ^{n+1}=\ZZ^\infty$ and  $\RK_1(C(\CC P^\infty_q))\simeq \{0\}$.
\end{thm}

\begin{proof}
Once again let us set $A_n = \K_0(C(\CC P^n_q))$ and $B_n = \K_1(C(\CC P^n_q))$, so that $A_n=\ZZ^{n+1}$ and $B_n=\{0\}$. 
Invoking the Milnor $\lim^1$-sequence we get 

$$0\map \lim^1_n B_n\map \RK_0(C(\CC P^\infty_q))\map \lim_n A_n\map 0.$$ From the argument of the above Proposition we conclude that the induced homomorphism $A_n\map A_{n-1}$ corresponds to the surjective map $\ZZ^{n+1}\map\ZZ^n$, i.e., projection onto the last $n$ summands. Consequently, $\lim_n A_n =\ilim_n \ZZ^{n+1}=\ZZ^\infty$ and clearly $\lim_n^1 B_n =\{0\}$. In the Milnor $\lim^1$-sequence for $\RK_1(C(\CC P_q^\infty))$, one finds $\lim_n B_n =\{0\}$ and the $\lim^1$-term vanishes, whence $\RK_1(C(\CC P^\infty_q))=\{0\}$.

\end{proof}

\noindent
Finally we turn our attention to the computation of the $\RK$-theory of $C(SU_q(\infty))$. It is known that $C(SU_q(n))$ is a type-$\textup{I}$ $C^*$-algebra for all $n\geqslant 2$ \cite{Br}, whence it is nuclear. The key step in our computation is the following result from \cite{N1}:

\begin{thm}[Nagy]
There are certain homological comparison elements $$\sigma_n\in\KK(C(SU(n)), C(SU_q(n))),$$ where $\KK$ denotes Kasparov's bivariant $\K$-theory, which induce isomorphisms $$\sigma_n: \K_i(C(SU(n))\overset{\sim}{\map}\K_i(C(SU_q(n))$$ for all $n$ and $i=0,1$. Moreover, there are commutative squares

\beq \label{naturality}
\xymatrix{
\K_{i}(C(SU(n)))\ar[r]^{\sigma_n} \ar[d] & \K_i(C(SU_q(n)))\ar[d]\\
\K_{i}(C(SU(n-1)))\ar[r]^{\sigma_{n-1}} & \K_i(C(SU_q(n-1)))
}
\eeq where the left (resp. right) vertical arrow is induced by the $*$-homomorphism $C(SU(n))\map C(SU(n-1))$ (resp. $C(SU_q(n))\map C(SU_q(n-1))$). 
\end{thm}

\begin{rem}
The homological comparison elements actually live in a bivariant $\K$-theory developed by Nagy, but this bivariant $\K$-theory agrees with Kasparov's $\KK$-theory for nuclear separable $C^*$-algebras. The above commutative diagram \eqref{naturality} follows from the explicit description of the homological comparison elements as partially defined $*$-homomorphisms on the generators $u^n_{ij}$'s (see Comment 3.9. of ibid.).
\end{rem}

\begin{prop} \label{KSUq}
$\RK_i (C(SU_q(\infty)))\simeq \ilim_n^1 \K_{1-i}(C(SU(n)))\oplus\ilim_n \K_{i}(CSU(n)))$.
\end{prop}

\begin{proof}
The Milnor $\lim^1$-sequence applied to the inverse system $\{C(SU_q(n))\}$ gives us the following short exact sequence

\beqn
0\map \lim^1_n \K_{1-i}(C(SU_q(n)))\map \RK_i(C(SU_q(\infty)))\map \lim_n \K_i(C(SU_q(n)))\map 0.
\eeqn From the above result of Nagy we conclude that $\ilim_n \K_i(C(SU_q(n)))\cong\ilim_n \K_i(C(SU(n)))$ and $\ilim^1_n \K_i(C(SU_q(n)))\cong\ilim^1_n \K_{i}(C(SU(n)))$. It is known that the $K$-theory of simply connected compact Lie groups is torsion free \cite{H}, whence $\lim_n \K_{i}(C(SU(n)))$ is torsion free. The above sequence splits and the assertion follows.
\end{proof}

\noindent
Let us set $\K^*(SU(n))=\K^0(SU(n))\oplus\K^1(SU(n))$. There is a $\ZZ/2$-graded Hopf algebra structure on $\K^*(SU(n))$, which is naturally isomorphic to $\K_*(C(SU(n)))$. Let $\rho: SU(n)\map U(N)$ be any unitary representation. Composing $\rho$ with the canonical inclusion $U(N)\hookrightarrow U(\infty)$ one obtains a map $\rho: SU(n)\map U(\infty)$, whose homotopy class determines an element of $\K^{-1}(SU(n))\cong\K^{1}(SU(n))$. Let $\rho^n_1,\cdots, \rho^n_{n-1}$ be the fundamental representations of $SU(n)$. Then we refer the readers to Theorem A of \cite{H}, a special case of which says

\begin{thm}[Hodgkin]
$\K^*(SU(n))\simeq \Lambda_\ZZ (\rho^n_1,\cdots, \rho^n_{n-1})$ as $\ZZ/2$-graded Hopf algebras.
\end{thm}

\noindent
The canonical inclusion $SU(n-1)\map SU(n)$ of Lie groups induces a $\ZZ/2$-graded homomorphism $\K^*(SU(n))\map\K^*(SU(n-1))$. The fundamental representations of $SU(n)$ admit a simple description, i.e., $\rho^n_i:SU(n)\map \Aut(\Lambda^i(V))$, where $V\simeq \CC^n$ denotes the standard representation of $SU(n)$. Using the branching rule of the restricted fundamental representation one finds that the induced $\ZZ/2$-graded ring homomorphism $\K^*(SU(n))\map \K^*(SU(n-1))$ sends 
\beqn
\text{$\rho^n_i\mapsto \rho^{n-1}_i +\rho^{n-1}_{i-1}$ for $1\leqslant i\leqslant n-1$ with $\rho^{n-1}_0=1$ and $\rho^{n-1}_{n-1} =1.$}\eeqn It follows that every generator $\rho^{n-1}_1,\cdots ,\rho^{n-1}_{n-2}$ of $\K^*(SU(n-1))$ has a preimage, whence the induced homomorphism is surjective. Consequently, the $\lim^1$-term in Proposition \ref{KSUq} above vanishes and one finds

\begin{thm}
$\RK_*(C(SU_q(\infty)))\cong \ilim_n \Lambda_\ZZ (\rho^n_1,\cdots , \rho^n_{n-1})$ as $\ZZ/2$-graded abelian groups. 
\end{thm}


\bibliographystyle{abbrv}

\begin{thebibliography}{10}


\bibitem{B}
B. Blackadar.
Shape theory for $C^*$-algebras.
{\em Math. Scand.} 56 (1985), no. 2, 249--275.

\bibitem{Bou}
N. Bourbaki.
General Topology, Chapters 1-4.
Elements of Mathematics, Springer-Verlag 1998.

\bibitem{Br}
K. Bragiel. 
The twisted $SU(N)$ group. On the $C^*$-algebra $C({S}_\mu{U}(N))$.
{\em Lett. Math. Phys.} 20 (1990), no. 3, 251--257.

\bibitem{CM}
A. L. Carey and J. Mickelsson. 
The universal gerbe, Dixmier-Douady class, and gauge theory.  
{\em Lett. Math. Phys.}  59  (2002),  no. 1, 47--60.

\bibitem{HM}
J. A. Harvey and G. Moore. 
Noncommutative tachyons and $\K$-theory. Strings, branes, and M-theory.  
{\em J. Math. Phys.}  42  (2001),  no. 7, 2765--2780.


\bibitem{H}
L. Hodgkin.
On the $K$-theory of Lie groups. 
{\em Topology} 6 (1967), 1--35.


%



\bibitem{KV}
J. Kustermans and S. Vaes. 
The operator algebra approach to quantum groups.  
{\em Proc. Natl. Acad. Sci. USA}  97  (2000),  no. 2, 547--552.

\bibitem{KVVVDW}
J. Kustermans, S. Vaes, L. Vainerman, A. Van Daele and S. Woronowicz. 
Lecture Notes School/Conference on Noncommutative Geometry and Quantum Groups, Warsaw 2001.
{\tt{https://perswww.kuleuven.be/~u0018768/artikels/lecture-notes.pdf}}


\bibitem{M}
A. Mallios.
Topological algebras. Selected topics.
North-Holland Mathematics Studies, 124, 1984.

\bibitem{N}
G. Nagy. On the Haar measure of the quantum $SU(N)$ group. Comm. Math. Phys. 153 (1993), 217--228.

\bibitem{N1}
G. Nagy. Deformation quantization and $K$-theory. Perspectives on quantization.
{\em Contemp. Math.} 214, Amer. Math. Soc., (1998), 111--134. 

\bibitem{P1}
N. C. Phillips. 
Inverse limits of $C^*$-algebras.  
{\em J. Operator Theory}  19  (1988),  no. 1, 159--195.

\bibitem{P2}
N. C. Phillips. 
Inverse limits of $C^*$-algebras and applications.  
Operator algebras and applications, Vol. 1, 127--185, London. Math. Soc. Lecture Note Ser., 135, 1988.

\bibitem{P3}
N. C. Phillips. 
Representable K-theory for $\sigma$-$C^*$-algebras.  K-Theory  3  (1989),  no. 5, 441--478. 

\bibitem{P4}
N. C. Phillips.
$K$-theory for Fr{\'e}chet algebras.
Internat. J. Math., Vol. 2 (1), (1991), 77--129.




\bibitem{Sheu}
A. Sheu.
Compact Quantum Groups and Groupoid $C^*$-Algebras
J. Func. Anal.
144,  (1997) 371-393.

\bibitem{Sheu2}
A. Sheu. 
Quantum spheres as groupoid $C^*$-algebras.  Quart. J. Math. Oxford Ser. (2)  48  (1997),  no. 192, 503Ð510.

\bibitem{Sheu3}
A. Sheu. 
Groupoid approach to quantum projective spaces.  Operator algebras and operator theory (Shanghai, 1997),  341--350, Contemp. Math., 228, Amer. Math. Soc., Providence, RI, 1998. 



\bibitem{S}
Ya. S. Soibelman.
Algebra of functions on a compact quantum group and its representations.
{\em Algebra i Analiz}, 2(1):190--212, 1990.

\bibitem{S-V}
Ya. S. Soibelman and L. L. Vaksman. The algebra of functions on the quantum group SU(n+1), and odd-dimensional quantum spheres. Leningrad Math. J. 2 (1991), 1023--1042.

\bibitem{Weid}
J. Weidner.  $KK$-groups for generalized operator algebras. I, II.  K-Theory  3  (1989),  no. 1, 57--77, 79--98.



\bibitem{W1}
S. L. Woronowicz.
A remark on compact matrix quantum groups.
{\em Lett. Math. Phys.}, 21 (1991), no. 1, 35--39.


\bibitem{W}
S. L. Woronowicz.
Compact quantum groups.
Sym{\'e}tries quantiques (Les Houches, 1995), 845--884, North-Holland, Amsterdam, 1998.














\end{thebibliography}

\vspace{5mm}
\noindent
\address{}

\end{document}